\documentclass[a4paper,10pt]{article}
\usepackage[utf8]{inputenc}
\usepackage[margin = 3cm]{geometry}
\usepackage{amsmath,amssymb,amscd,amsthm}
\usepackage{multicol}
\usepackage{latexsym}
\usepackage{mathrsfs}
\usepackage{fancyhdr}
\usepackage[dvips]{color}
\newcommand{\BigO}[1]{\ensuremath{\operatorname{O}\bigl(#1\bigr)}}
\newtheorem{mytheorem}{Theorem}
\newcommand{\Cliff}{\mbox{Cliff}}

\title{Counterexamples to Mercat's Conjecture}
\author{Akash Kumar Sengupta}

\begin{document}

\maketitle

\begin{abstract}
In this paper, we provide counterexamples to Mercat's conjecture on vector bundles of rank $n \geq 4$ on algebraic curves. 
For any $n\geq 4$, we provide examples of curves lying on K3 surfaces and vector bundles of rank $n$ on those
curves for which Mercat's inequality fails.\\
\end{abstract}

\section{Introduction}
The Clifford index $\gamma_C$ (or Cliff($C$)) of a curve $C$ measures the complexity of the curve in its moduli space. Lange and Newstead (see [2]) proposed the following definition as a generalization
of $\gamma_C$ for higher rank vector bundles. If $E \in \mathscr{U}_C(n,d)$ is a semi-stable vector bundle of 
rank $n$ and degree $d$ on a curve $C$ of genus $g$, its Clifford index is defined as :
\[ \gamma(E) := \mu(E) - \frac {2} {n} h^0(C,E) + 2 \geq 0 ,\]
and the higher Clifford indices of $C$ are defined as:
\[ \Cliff_n(C) = min \{ \gamma(E) : E \in \mathscr{U}_C(n,d), d < n(g-1), h^0(C,E) \geq 2n \}. \]
We have, Cliff$_1(C) = \gamma_C$. Considering direct sums of line bundles, it is easy to see that Cliff$_n(C) \leq \gamma_C$ for all $n\geq 1$ (see [11]). Mercat
(see [3]) proposed a conjecture relating these Clifford indices to the classical geometry of $C$. It was rephrased in the following form in [2] (Conj. 9.3):
\[(M_n) : \quad \Cliff_n(C) = \gamma_C .\]
Mercat's conjecture $(M_2)$ is known to be true for various classes of curves (see [2]), in particular for general
$k$-gonal curves of genus $g > 4k - 4$. Farkas and 
Ortega provided counter examples to $(M_2)$ in [11]. They showed that for every genus $g\geq 11$ there exist curves $[C] \in \mathscr{M}_g$ 
with maximal Clifford index carrying stable rank 2 bundles $E$ with 4 sections such that $\gamma(E) < \gamma_C$. 
\paragraph{}
Lange, Mercat and Newstead provided counterexamples to $(M_3)$ by constructing rank 3 bundles with 6 independent
sections on a general curve of genus 9 or 11 (see [14]). Also, it was shown in [5] that for sufficiently high genus
Mercat's conjecture $(M_3)$ fails for any smooth curve of maximal Clifford index lying on a K3 surface.
\paragraph{}
Counterexamples to Mercat's conjecture $(M_4)$ were provided in [6]. The authors constructed rank 4 
bundles on curves lying on K3 surfaces of picard number 1 for which Mercat's inequality fails.
\paragraph{}
In this paper, we provide counterexamples to $(M_n)$ for all $n\geq 4$. Our main result is the following:
 \begin{mytheorem} For each integer $n\geq4$, there exists a $\mathrm{K3}$ surface $S$ with a stable vector bundle $E$ on $S$ and a smooth curve $C\subset S$ such that,
\begin{itemize}
 \item[(0)] $E|_C$ is a semi-stable vector bundle on $C$,
 \item[(1)]  $ \mathrm{rk}(E|_C) = n$,
 \item[(2)] $ \mathrm{deg}(E|_C) < n(g(C) -1)$,
 \item[(3)] $h^0(E|_C) \geq 2n$,
 \item[ and,]
 \item[(4)] $\mathrm{Cliff}_n(C) \leq \gamma(E|_C) < \gamma_C$.
\end{itemize}
Hence Mercat's conjecture $(M_n)$ fails for $C$.
 \end{mytheorem}
\paragraph{}For a K3 surface $S$, a smooth curve $C \subset S$ and a globally generated linear series $A \in \displaystyle W_d^r(C)$ with $h^0(C,A)=r+1$,
the Lazarsfeld-Mukai bundle $E := E_{C,A}$ is defined via the following elementary modification on $S$
\[ 0 \longrightarrow E_{C,A}^\vee \longrightarrow H^0(C,A)\otimes \mathscr{O}_S \longrightarrow A \longrightarrow 0 . \]
In [6], the authors related Mercat's conjecture $(M_{r+1})$ to semistability of the restricted Lazarsfeld-Mukai bundle $E|_C$ on the curve $C$
and proved semistability in rank 4 to obtain counterexamples to $(M_4)$. However the techniques used in [6] and other previously
mentioned counterexamples were rank-specific and did not seem to generalize to higher ranks. In the proof of Theorem 1, corresponding to each integer 
$n\geq 4$, we choose a triple $(S,C_0,A)$ suitably such that $\mathrm{Pic}(S)=\mathbb{Z}.[C_0]$ and the LM bundle $E =E_{C_0,A}$ is $C_0$-stable.
Using Flenner's restriction theorem, we choose a curve $C$ (different from $C_0$) lying on $S$ such that the restricted bundle
$E|_C$ is semi-stable. Then it turns out that the bundle $E|_C$ contributes to the n-th Clifford index of $C$ and indeed $\gamma(E) < \gamma_C$. Thus the bundle $E|_C$ on the curve $C$ invalidates
Mercat's conjecture in rank $n$.
\paragraph{}
The structure of the paper is as follows; in section 2, we recall basic properties of Lazarsfeld-Mukai bundles on K3 surfaces. Finally, in section 3, we prove our main result providing counterexamples
 to Mercat's conjecture.\\\\
 {\bf Acknowledgements. } I am grateful to Gavril Farkas and Angela Ortega for fruitful discussions and encouragement. I am thankful to
Peter Newstead for comments on the preliminary version of the paper.
I would like to thank Berlin Mathematical School and Humboldt Universit\"{a}t zu Berlin for the hospitality during the preparation of this work.

\section{Lazarsfeld-Mukai Bundles}
Let $S$ be a smooth, complex projective  K3 surface and $L$ a globally generated line bundle on $S$ with $L^2= 2g-2$. Let $C\in |L|$ be a smooth curve and $A \in W_d^r(C)$ a base
point free line bundle with $h^0(C,A) = r+1$. Considering $A$ as a globally generated sheaf on S, define the bundle $F_{C,A}$ to be the kernel of the evaluation map 
$ev_{A,S} : H(C,A)^0 \otimes \mathscr{O}_S \mapsto A$, i.e.,
\[ 0 \longrightarrow F_{C,A} \longrightarrow H^0(C,A)\otimes \mathscr{O}_S \longrightarrow A \longrightarrow 0 .\]
The Lazarsfeld-Mukai bundle associated to the pair $(C,A)$ is, by definition, $E_{C,A} := F_{C,A}^\vee$. We recall the following properties of the bundle $E_{C,A}$ below (see [7], [8]):
\begin{itemize}
 \item[(i)] $\mathrm{rk}(E_{C,A}) = r+1$,
\item[(ii)] $\mathrm{det}(E_{C,A})= L$,
\item[(iii)] $c_2(E_{C,A}) = d$,
\item[(iv)] $E_{C,A}$ is globally generated off the base locus of $\omega_C \otimes A^\vee$,
\item[(v)] $h^0(S,E_{C,A}) = h^0(C,A) + h^0(C,\omega_C \otimes A^\vee) = r+1+g-d+r$,
\\ and $h^1(S,E_{C,A}) = h^2(S,E_{C,A}) = 0$,
\item[(vi)] $\chi(S,E_{C,A} \otimes F_{C,A}) = 2(1-\rho(g,r,d)$.
\end{itemize}

When $\mathrm{Pic}(S) = \mathbb{Z}.L$ and $\rho(g,r,d) =0$, the it is well-known that there is only one single Lazarsfeld-Mukai bundle $E$ with $c_1(E) = L$ and $c_2(E) =d$.
Note that in this case $\rho(g,r',d') = 0$, where $r'= h^0(\omega_C \otimes A^\vee) -1$, $d' = deg(\omega_C \otimes A^\vee)$. Hence $\omega_C \otimes A^\vee$ is base point free
and $E_{C,A}$ globally generated.
\paragraph{ Lemma 2.1: } If $\mathrm{Pic}(S) = \mathbb{Z}.L$ and $\rho(g,r,d) =0$, then the bundle $E:= E_{C,A}$ is $\mu$-stable for the polarization given by $L$ (or $C$-stable, in short).
\begin{proof} Suppose $E$ is not stable. Have a destabilizing sequence
\[ 0 \longrightarrow M \longrightarrow E \longrightarrow N \longrightarrow 0 ,\]
where $M$ and $N$ are sheaves on $S$ with first Chern classes $c_1(M)=mL$ , $c_1(N)=nL$ for some $m,n \in \mathbb{Z}$. 
As noted above, $E$ is globally generated and hence so is $N$. Therefore $n\geq 0$.\\

If $n=0$, then the torsion-free part of $N$ is a trivial bundle. As $\mathrm{rk}(M) < \mathrm{rk}(E)$, the torsion-free part 
of $N$ is non-zero. Hence $E$ has a trivial quotient, say $N_1$. Then $N_1^\vee$ has global sections and consequently 
$E^\vee$ has global sections. But, via Serre duality, $\displaystyle h^0(S,E^\vee) = h^2(S,E) =0$. Therefore $n>0$.\\

Now $L = c_1(E) = c_1(M) + c_1(N) = (m+n)L$, hence $m = 1 -n \leq 0$. Then we have
\[\mu(M) = \frac{mL^2}{\mathrm{rk}(M)} \leq 0 < \frac{L^2}{\mathrm{rk}(E)} = \mu(E) ,\]
which contradicts the assumption of $M$ being destabilizing.
\end{proof}

\section{Proof of Theorem 1}
%As described in the introduction, we have the following setup.
Let $r\geq 3$ , we set $n:= r+1$ and $g :=(r+1)(r+2)$.\\\\
Let $S$ be a K3 surface with a line bundle $L$ such that,
\begin{itemize}
\item[(i)] $Pic(S)= \mathbb{Z}.L$,
\item[(ii)] $L$ is ample and globally generated,
\item[(iii)] $|L|$ has a smooth member and $L^2 = 2g-2$
\end{itemize}
As $g>2$, there exists such a pair $(S,L)$ by standard arguments (see [12]).
\paragraph{}
Let $C_0\in |L|$ be a general smooth curve. It follows from a theorem of Lazarsfeld (see [7]) that $C_0$ is Brill-Noether generic.
So $\rho(g,r,d) = g - (r+1)(g-d+r) \geq 0$ if and only if the Brill-Noether locus $W_d^r(C_0)$ is non-empty (see [1]).
\paragraph{}
Let $\displaystyle d := g-[\frac{g}{r+1}] +r = (r+1)(r+2)-2$, i.e. $d$ is minimal such that $W_d^r(C_0)$ is non-empty.
Let $A \in W_d^r(C_0)$ be a linear series.
As $d$ is minimal, we have $h^0(C_0,A)=r+1$ and $A$ is globally generated ([2], Lemma 4.2).
\paragraph{}Let $E := E_{C_0,A}$ be the Lazarsfeld-Mukai bundle corresponding to the pair $(C_0,A)$. 
Using the results in section 2, we have,
\begin{itemize}
 \item[(E1)] $\mathrm{rk}(E) = n$,
\item[(E2)] $\mathrm{det}(E)=L$,
\item[(E3)] $c_2(E) =d$,
\item[(E4)] $h^0(E) = r+1+g-d+r = r+1+r+2 > 2r+2 = 2n$,
\item[(E5)] as $\mathrm{Pic}(S) = \mathbb{Z}.L$, we have $E$ is $C_0$-stable.
\end{itemize}
\paragraph{}
Let $a := (n^2-1)(g-1)$, and $C\in|aL|$ be a general smooth curve.
Then
\[\displaystyle\frac{{ a+2 \choose 2} - a -1}{a} > (2g-2)(\frac{n^2-1}{4}),\]
and hence by Flenner's restriction theorem ([10], Theorem 7.1.1), $\displaystyle E|_C$ is semi-stable vector bundle on $C$ with $rk(E|_C) = n$. So we have parts (0) and (1) of the Theorem.
\paragraph{}Now, by the adjunction formula, $C$ is of genus
\[\displaystyle g(C) = 1+ \frac{a^2L^2}{2}=1+a^2(g-1).\]
Then we have
\[\mathrm{deg}(E|_C) = ac_1(E).L = aL^2 = a(2g-2) < na^2(g-1) = n(g(C) - 1),\]
 which proves assertion (2) of the Theorem.\\\\
As noted in (E4), we have $h^0(E) > 2n$. Hence, part (3) of the theorem will follow from the next lemma.\\\\
{\bf Lemma 3.1} : $h^0(E) \leq h^0(E|_C)$.
\begin{proof} Have the short exact sequence,
\[ 0\longrightarrow \mathscr{O}_S(-C) \longrightarrow \mathscr{O}_S \longrightarrow \mathscr{O}_C \longrightarrow 0.\]
Tensoring by $E$ and taking cohomologies we obtain,
\[ 0 \longrightarrow H^0(E(-C)) \longrightarrow H^0(E) \longrightarrow H^0(E|_C) \longrightarrow \dots\]
\paragraph{}If $H^0(E(-C)) \neq 0$, then $\mathscr{O}_S(C) \hookrightarrow E$.
Since $E$ is $C_0$-stable, we have
\[\mu(\mathscr{O}_S(C)) = aL^2 < \mu(E) = \frac{L^2}{n},\]
which is a contradiction! Hence $H^0(E(-C)) = 0$ and $h^0(E) \leq h^0(E|_C)$.\\
\end{proof} 
%\paragraph{}Therefore $h^0(E|_C) \geq h^0(E) > 2n$, and hence we have (3).
\paragraph{} Parts (0)-(3) of the theorem imply that the semistable bundle $E|_C$ contributes to the $n$-th Clifford index
$\mathrm{Cliff}_n(C)$ of the curve $C$. Now suppose Mercat's conjecture $(M_n)$ is true. Then we have,
\[\displaystyle\gamma_C = \Cliff_n(C) \leq \mu(E|_C) - \frac{2}{n} h^0(E|_C) + 2 \]
\hspace*{4cm}$\displaystyle\Rightarrow \gamma_C -2 + \frac{2}{n} h^0(E|_C) \leq \mu(E|_C)$\\
\newline
\hspace*{4cm}$\displaystyle\Rightarrow \gamma_C +2 \leq \mu(E|_C) = \frac {aL^2}{n}$, as $h^0(E|_C)> 2n$.\\
\newline
Hence, $\displaystyle(M_n) \Rightarrow \gamma_C +2 \leq \frac {aL^2}{n}$. \hspace*{6.7cm} $(*)$
\paragraph{}Next we compute the Clifford index of $C$.
Consider the line bundle $L|_C$ on $C$.
Have the short exact sequence:
\[ 0\longrightarrow L \otimes \mathscr{O}_S(-C) \longrightarrow L \longrightarrow L \otimes \mathscr{O}_C \longrightarrow 0 .\hspace*{2cm}(**)\]
As $C \in |aL|$, we have $H^0(L \otimes \mathscr{O}_S(-C)) = 0$.
Then, taking cohomologies in $(**)$ and by using Riemann-Roch we get 
\[\displaystyle h^0(L|_C) \geq h^0(L) = \frac {L^2} {2} + 2 \geq 2,\]
and $\mathrm{deg}(L|_C) \leq g(C) -1$. Hence $L|_C$ contributes to 
the Clifford index of $C$ and
\[ \gamma_C \leq \mathrm{deg}(L|_C) - 2(h^0(L|C) - 1) \leq (a-1)L^2 -2 < \lfloor \frac{g(C) -1} {2}\rfloor.\]
%We have the following two cases depending on the Clifford index of $C$:\\\\
%{\bf Case I:} $\gamma_C = \lfloor \frac{g(C) -1} {2}\rfloor$.\\\\
%Then, $\gamma_C = \lfloor \frac{a^2L^2}{4} \rfloor \geq  \frac{a^2L^2}{4} -1$ and by (*), $\gamma_C +2 \leq \frac {aL^2}{n}$.\\\\
%Hence, $\frac{a^2L^2}{4} +1 \leq \frac {aL^2}{n}$, which is a contradiction since $n\geq 4$ and $a\geq 2$.\\\\
%{\bf Case II:} $\gamma_C < \lfloor \frac{g(C) -1} {2}\rfloor$.\\\\
\paragraph{}Hence $C$ is not of maximal Clifford index.
Then it follows from the main theorem of Green-Lazarsfeld ([4]), that there is a line bundle $M$ on $S$ whose restriction to $C$ computes the Clifford
index of $C$ and we have the following (see [9], Lemma 8.3):
\[\gamma_C = deg(M \otimes \mathscr{O}_C) - 2(h^0(M) -1) = M.\mathscr{O}_S(C) - M^2 -2.\]
Since $\mathrm{Pic}(S) =\mathbb{Z}.L$, we have $M = xL$ for some $x\in \mathbb{Z}$.
\paragraph{}Then $\gamma_C = xL.aL - x^2L^2 -2 = x(a-x)L^2 -2$. Hence the minimum is attained when $x=1$, i.e. the Clifford index of $C$ is computed by $L$ itself and 
$\gamma_C = (a-1)L^2 -2$.
Then by $(*$),
\[\gamma_C + 2 \leq \frac {aL^2}{n} \Rightarrow (a-1)L^2 \leq \frac {aL^2}{n} \Rightarrow n(a-1) \leq a,\]
which gives a contradiction since $n\geq 4$ and $a\geq 2$. 
\paragraph{}This completes our proof showing that Cliff$_n(C) \leq \gamma(E|_C) < \gamma_C$ and Mercat's conjecture $(M_n)$ fails for $C$. 
\paragraph{ Remark 3.2:} The vector bundle
$E|_C$ has Clifford index $\displaystyle\gamma(E|_C) \leq \mu(E|_C) -2 = \frac{a(2g-2)}{n} -2$ and the Clifford index of the curve is $\gamma_C \geq (a-1)(2g-2)$.\\\\
Therefore, $\displaystyle\gamma_C - \Cliff_n(C) \geq (a-1)(2g-2) - \frac{a(2g-2)}{n} +2 = \BigO {n^6}$.
\paragraph{}However, the discrepancy is not as huge as the genus $g(C) = \BigO{n^{10}}$. It might be interesting to find counter-examples with maximum
possible difference $\gamma_C - \Cliff_n(C)$, as done in the case of rank 2 in [13].
\paragraph{Remark 3.3:} Theorem 1 shows that Mercat's conjecture fails in rank $n$ for the curve $C$, but the curve here is not a general curve. 
The conjecture $(M_n)$ on a general curve of arbitrary genus $g$ is still unanswered in most cases, even in the case of rank 2 bundles.

\paragraph{}{\bf References}
\begin{itemize}
\item[{[1]}] E. Arbarello, M. Cornalba, P.A. Griffiths, J. Harris, \emph {Geometry of algebraic curves} 
Vol I Springer-Verlag, New York, 1985.
\item[{[2]}] H. Lange, P.E. Newstead, \emph{Clifford indices for vector bundles on curves}
, arXiv:0811.4680, in: Affine
Flag Manifolds and Principal Bundles (A. H. W. Schmitt editor), Trends in Mathematics, 165-202,Birkhauser (2010).
\item[{[3]}] V. Mercat, \emph{Clifford’s theorem and higher rank vector bundles}, International Journal of Mathematics 13 (2002), 785-796.
\item[{[4]}] M. Green, R. Lazarsfeld, \emph{Special divisors on curves on a K3 surface}, Inventiones Math.89 (1987),357-370.
\item[{[5]}] G. Farkas, A. Ortega, \emph {Higher Rank Brill-Noether theory on sections of K3 surfaces},International Journal of Mathematics
Vol. 23, No. 7 (2012) 1250075. 
\item[{[6]}] M. Aprodu, G. Farkas, A. Ortega, \emph {Minimal resolutions, Chow forms of K3 surfaces and Ulrich bundles},  arXiv:1212.6248v3 18 Apr 2013 
\item[{[7]}] R. Lazarsfeld \emph{Brill-Noether-Petri without degenerations},  J. Differential Geom. Volume 23, Number 3 (1986), 299-307. 
\item [{[8]}] M. Aprodu, \emph {Lazarsfeld-Mukai bundles and applications}, arXiv:1205.4415v1 [math.AG]. 20 May 2012
\item [{[9]}] A. L. Knutsen, \emph{On kth-order embeddings of K3 surfaces and Enriques surfaces}, Manuscripta Mathematica, February 2001, Volume 104, Issue 2, pp 211-237, (arXiv:math/0008169v2 [math.AG])
\item [{[10]}] D. Huybrechts and M. Lehn, \emph{The geometry of the moduli space of sheaves}, Second Edition, Cambridge University Press 2010
\item [{[11]}] G. Farkas, A. Ortega, \emph{The maximal rank conjecture and rank two Brill-Noether theory}, Pure and Applied Mathematics Quarterly 7 (2011), 1265-1296  volume dedicated to Eckart Viehweg ( ArXiv: math.AG/ 1010.4060, preprint).
\item [{[12]}] A. Beauville, \emph{Complex Algebraic Surfaces}, London Mathematical Society Student Texts 34 
\item [{[13]}] H. Lange, P. E. Newstead, \emph{Bundles of rank 2 with small Clifford index on algebraic curves},  arXiv:1105.4367v2 [math.AG], 22 May 2011
\item [{[14]}] H. Lange, V. Mercat, P. E. Newstead, \emph{On an example of Mukai}, Glasgow Math J. 54 (2012) 261-271\\
\end{itemize}
\hspace*{1cm}CHENNAI MATHEMATICAL INSTITUTE, H1 SIPCOT IT PARK, \\
\hspace*{1cm}CHENNAI, INDIA\\
\hspace*{1cm}E-mail address: $asengupta@cmi.ac.in$
\end{document}